Mathematics Hidden Behind the Two Coefficients of Babylonian Geometry

Kazuo MUROI

§1. Introduction

In Babylonian mathematics many geometrical coefficients are listed in some cuneiform tablets which we call coefficients lists. Some of these may serve to evaluate the level of Babylonian mathematics if we could see through the mathematical meanings of the numerical values and their modifiers listed together. For example, the following two items would be representative of the geometrical coefficients:

1,24,51,10 *ṣi-li-ip-tum* íb-$si_8$ "1;24,51,10 is the diagonal (of) a square"

4,48 *ku-bu-ur i-ṣí-im* "4,48 is the thickness of a log".[1]

As Neugebauer and Sachs clarified, the former means:

√2 ≈ 1;24,51,10 (=1.414212962 ⋯)

and the latter substantially means:

1/4π ≈ 0;4,48 (that is, π ≈ 3.125),

both of which show that the Babylonian scribes certainly had the ability to cope with rather complicated calculations.

In lines 5 and 6 of Susa mathematical text no. 3,[2] one of the Old Babylonian coefficients lists, we have two enigmatic items which seem to concern the areas of certain circular figures:

16 [igi]-gub *šà* gúr *šà* 2 še *i-na* šà gúr gar

"0;16[sic] is the coefficient of a circular figure in the midst of which two barley (figures) are put." (line 5)

16,26,46,40 igi-gub *šà* gúr *šà* 3 še *i-na* šà gúr gar

"0;16,26,46,40 is the coefficient of a circular figure in the midst of which three

barley (figures) are put." (line 6)

At a glance we may notice that the second numerical value of four figures, which is rare in coefficients lists, involves a complex geometrical figure or complicated calculations. However, the mathematical meaning of line 6 as well as that of line 5 has wholly been unclear since the text was published in 1961. In the present paper I shall reveal two circular figures hidden behind the above-mentioned two items of the Susa coefficients list with my own analysis of the text.

§2. Technical terms

Before making a mathematical analysis of the two lines in question, we should clarify the literal and mathematical meanings of the three technical terms which occur in the text.

(1) igi-gub "coefficient"

The literal meaning of igi-gub or igi-gub-ba is probably "that which stands by (you in your calculations)"[3]. This term occurs in the so-called problem texts as well as the coefficients lists.

(2) gúr "circle, circumference"

The Sumerian noun gúr "loop, hoop, circle" is derived from the verb gúr "to bend, curve" and its literal meaning would be "the curved". In mathematical texts it usually means "circle", but the translation of "circle" does not fit the numerical values of our two coefficients. I propose to translate our gúr as "circular figure", which will be confirmed by a mathematical analysis below.

(3) še "barley"

The technical term še "barley (figure)", which also occurs in lines 16-18 of Susa mathematical text no. 3, is a symmetrical figure surrounded by two quadrants. In two other coefficients lists and one problem text the same figure is also called gán-(giš)má-gur$_8$ "field of a deep-going boat".[4]  See fig. 1.

§3. Analysis of line 5

After many attempts to interpret the difficult expression that modifies the number 16[sic] in line 5 by the usual translation "circle" of gúr, I have reached the conclusion that the figure in question is not a circle but a circular figure composed of two "barley figures" which overlap each other. See fig. 2. This is a famous figure, at least Japanese junior high schools, because the calculation of the area of this figure is not easy for almost all the students, though it seems to be elementary. However it is within the reach of the Babylonian scribes, the fact I know from my own experience studying mathematical texts for a long time. Since the area of the circular figure in line 5 (= S) is $\pi/3 + 1 - \sqrt{3}$ and the approximate values of $\pi$ and $\sqrt{3}$ are usually 3 and 7/4 (= 1;45) respectively in Babylonian mathematics, it follows that S ≈ 0;15. The number 16 in line 5 may be a scribal mistake for 15 which was possibly caused by the occurrence of the same number at the head of line 6.[5]

Now, let us try to reconstruct the Babylonian method for calculating S, assuming the radius of a quadrant to be 1. See fig. 3.

Since ∠ACB= 30°, AD = 1/2 and BD = 1 − √3/2, we obtain

$AB^2 = AD^2 + BD^2 = 2 - \sqrt{3}$ ($\approx 2 - 1;45 = 0;15$), which is the area of the square inscribed in the circular figure.

Also, the area of the segment ABE is:

$(1/12)T - (1/2) \cdot 1 \cdot (1/2) = T/12 - 1/4$, where T is the area of the circle whose radius is 1.

Therefore, $S = 2 - \sqrt{3} + 4(T/12 - 1/4) = T/3 + 1 - \sqrt{3}$

$\approx 3/3 + 1 - 1;45 = 0;15$.

It should be noted that S equals $AB^2$ in the end because of $T \approx 0;5 \cdot 6^2 = 3$.

§4. Analysis of line 6

The circular figure designated in line 6 as "that in the midst of which three barley (figures) are put" is a regular hexagonlike figure composed of six equilateral triangles and six small circular segments. See fig. 4. In calculation of the area of this figure the scribe of our tablet seems to have assumed the radius of a quadrant to be 1 again and have omitted the six small segments, because one sixth of the numerical value listed in line 6 is:

$(1/6) \cdot 0;16,26,46,40 = 0;2,44,27,46,40 = 0.0\mathbf{456}84\cdots$, which is close to the area of an equilateral triangle whose side is $(\sqrt{14} - \sqrt{6})/4$ (= s):

$(\sqrt{3}/4)s^2 \approx (7/16)s^2 = 0.0\mathbf{456}55\cdots$,

and the area of the segment is approximately 0.0028, which is larger than the difference between the two values above.

In this way, the scribe approximated the area in question to the total area of the six

equilateral triangles. The length of the side s is obtained as follows (see fig. 5):

Since BH = s/2, HO = √3s/2, OC = √2/2, and CB = 1, we obtain the quadratic equation $(s/2)^2 + (\sqrt{3}s/2 + \sqrt{2}/2)^2 = 1$, that is, $s^2 + (\sqrt{6}/2)s = 1/2$.

By completing the square, which was called *takīltum* in Akkadian, the solution of the equation is accomplished;

$(s + \sqrt{6}/4)^2 = 14/16$, and so $s = (\sqrt{14} − \sqrt{6})/4$.

Therefore, according to the Babylonian formula, the area of the equilateral triangle is:

$(7/16)s^2 = 7(5 − \sqrt{21})/64$.

As we have seen above, this value (0.**0456**55⋯) is a little smaller than one sixth of the value in line 6 (0.**0456**84⋯) and this fact may suggest that the approximate value of √21 used by the scribe is a little smaller than the correct value √21 = 4;34,57,16,21,0,23, ⋯. If we suppose, for example, √21 to be 4;34,57,15,10,28, we obtain the following value:

$7(5 − \sqrt{21})/64 \approx 0;\mathbf{2,44,27,46,40},1,25$

$= 0.045684\cdots$.

It seems that the scribe made some mistake in the routine though complicated calculation of the square root of 21. It would be instructive to see how the approximation of √21 is obtained by the Babylonian method:

Since $4 < \sqrt{21} < 5$, the first approximation is $a_1 = (4 + 5)/2 = 4;30 < \sqrt{21}$, and so the second is $a_2 = 21/a_1 = 14/3 = 4;40 > \sqrt{21}$.

Applying the same method to $a_1$ and $a_2$, we obtain $a_3 = (a_1 + a_2)/2 = 55/12 = 4;35 > \sqrt{21}$, and $a_4 = 21/a_3 = 252/55 < \sqrt{21}$.

Moreover, $252/55 = 4;34,54,32,43,\cdots$.

Therefore, we immediately obtain the fifth approximation, which is correct to five figures:

$a_5 = (a_3 + a_4)/2 = 4;34,57,16,21,\cdots$

The scribe must have made a mistake in the final step of this calculation or in dividing 252 by 55.

§5. Conclusion

Through making a mathematical and linguistic analysis of the two enigmatic items of the Susa coefficients list and clarifying the two circular figures hidden behind them, we have reconfirmed that the Babylonians left behind great achievements in geometry as well as in algebra, although their main interests were the calculations of the areas of various figures. In addition, we have realized, I believe, the scribe of this Susa tablet clearly excels in mathematics average modern people, of course including me. I would like to respect him for his good sense of geometrical figures.

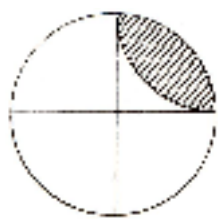

Fig. 1  še or gúr-má-gur₄

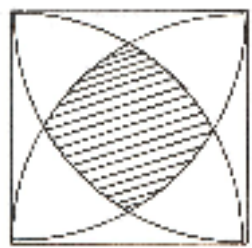

Fig. 2. gúr ša 2 še 1-na šà gúr gor

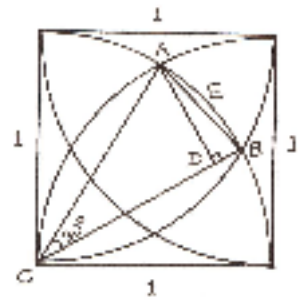

Fig. 3.

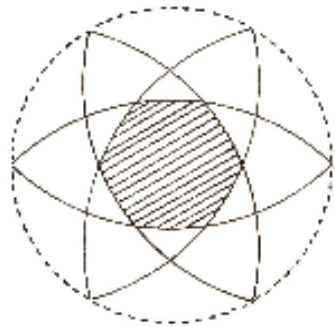

Fig. 4. gúr šà 3 še i-na šà gúr ŋar

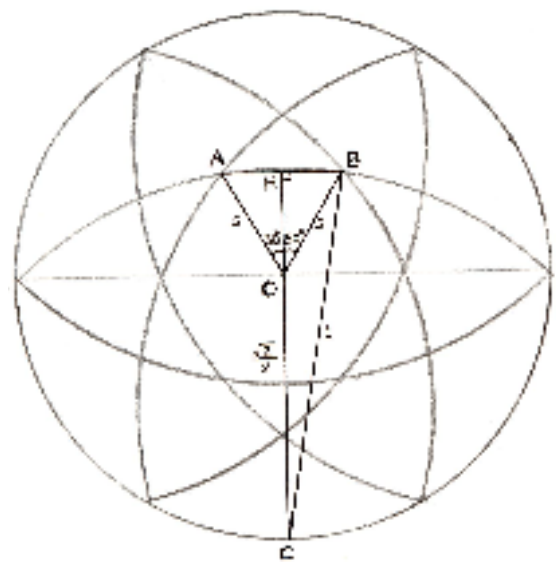

Fig. 3.

Notes

(1) O. Neugebauer and A. Sachs, *Mathematical Cuneiform Texts* [= MCT], 1945, pp. 136-139, YBC 7243, lines 10 and 35.

(2) E. M. Bruins et M. Rutten, *Textes mathématiques de Suse* [= TMS], pp. 25-34.

(3) For an example of the compound verb igi ⋯ gub "to serve, to be at the service of", see:

M. T. Roth and others, *The Assyrian Dictionary* (= CAD), vol. 20, U and W, 2010, p. 382, 7a.

(4) CAD, M part 1, 1977, pp. 141-142.

(5) It is regrettable that there is no photograph of this tablet in TMS. I cannot help believing that Rutten's hand-copy of the tablet is correct.